\documentclass[preprint,12pt]{elsarticle}

\usepackage{amssymb}
 \usepackage{amsthm}
 \usepackage{latexsym}
\usepackage{amsmath}

\newtheorem{Theorem}{Theorem}[section]
\newtheorem{corollary}[Theorem]{Corollary}
\newtheorem{Definition}[Theorem]{Definition}

\newcommand{\Imp}{\mbox{$\Longrightarrow$}}

\def\ie{\emph{i.e.}}
\def\pes{\emph{e.g.}}
\def\Pes{\emph{E.g.}}

\def\F{{\mathcal F}}

\def\P{{\mathcal P}}

\def\U{{\mathcal U}}

\def\V{{\mathcal V}}
\def\W{{\mathcal W}}

\def\C{{\mathcal C}}

\def\N{{\mathbb N}}

\def\R{{\,\mathbb R}}

\def\hR{{}^{*}\R }

\def\hM{{}^{*}\! M }
\def\hP{{}^{*}\! P }
\def\hRR{{}^{*}\! R }

\def\g{{}^{*}\! g }

 \def\p{{}^{*}\! p}
 \def\q{{}^{*}\! q}
 \def\f{{}^{*}\! f}
 
 \def\X{{}^{*}\! X}
 \def\hA{{}^{*}\! A}
 \def\ns{{}^{*}\!}

	\def\B{\mathcal{B}}

\def\bg{{\beta}}

\def\yg{{\upsilon}}

\def\+#1{\vec{#1}}

\def\as{{\mathsf a}}
\def\cs{{\mathsf c}}

\def\es{{\mathsf e}}

\def\is{{\mathsf i}}

\def\ps{{\mathsf p}}

\def\uf{{\mathfrak{u}}}

\def\*{\times}
\def\st{{}^{\star}\!\, }
\def\0{\emptyset}
\def\/{\setminus}
\def\_{\overline}

\def\<{\prec}

\def\incl{\subseteq}

\def\transfer{\textit{Transfer Principle}}
\def\zfc{\textsf{ZFC}}

\def\qed{${}$\hfill $\Box$}

\def\trans{transfer principle}

\def\tran{\textbf{Tran}}
\def\ind{\textbf{Ind}}
\def\poss{\textbf{Poss}}

\journal{Topology and its Applications}

\begin{document}

\begin{frontmatter}



\title{Leibniz's Principles and Topological Extensions}


\author{Marco Forti}

\address{Dipart. di Matematica Applicata ``U. Dini'',
Universit\`{a} di Pisa, Italy. \tt{forti@dma.unipi.it}}

\begin{abstract}
    
    Three philosophical principles are often quoted in connection with Leibniz:
    ``objects sharing the same properties are the same object'', 
     ``everything can        
       possibly exist, unless it yields contradiction'',
     ``the ideal elements          
      correctly determine the real things''.

      Here we give a precise formulation 
       of these principles within the framework of
      the Topological Extensions of \cite{DF05}, structures     
       that generalize at once compactifications, completions, and nonstandard          %
       extensions. In this topological context, the above Leibniz's principles 
	appear  
       as a property of separation, a property of compactness, and 
       a property of analyticity, respectively. 
      
   Abiding by this interpretation, we obtain the somehow surprising 
   conclusion that these Leibnz's principles 
    \emph{can be fulfilled in pairs, but not all          
       three together}.

\end{abstract}


\end{frontmatter}


\section*{Introduction}

Three philosophical principles are often quoted in connection with Leibniz:

      \textbf{Identity of indiscernibles}
      \begin{center}
          \emph{``objects sharing the same properties are the 
               same object''}
      \end{center} 
 
     \begin{quote}
    
    There are never in nature two beings which are perfectly identical to 
    each other, and in which it is impossibile to find any internal difference
     \ldots\ (\emph{Monadology})

Two  indiscernible individuals cannot exist. [\ldots]
To put two 
indiscernible things is to put the same thing under two names.

(\emph{Fourth letter to Clarke})

[...] dans les choses sensibles on n'en trouve 
jamais deux indiscernables [...]\footnote{~among sensible things one 
never finds two [that are] indiscernible.}

(\emph{Fifth letter to Clarke}, \cite{CL}, p. 132)

\end{quote}

\textbf{Possibility as consistency}

     \begin{center}
          \emph{``everything can        
            possibly exist, unless it yields contradiction''}
     \end{center}
   
   \begin{quote}
    
     - \emph{Impossible} is what yields an absurdity.
    
     - \emph{Possible} is not impossible.
    
     - \emph{Necessary} is that, whose opposite is impossible.
    
     - \emph{Contingent} is what  is not  necessary.
     
    (\emph{unpublished}, 1680 ca.)

    \bigskip
    
    [\ldots] nothing is absolutely necessary, when the contrary 
    is possible.\ \ \ [\ldots]\\
    Absolutely necessary is [\ldots] that whose opposite yields a 
contradiction. 

(\emph{Dialogue between Theofile and Polydore}) 

    \end{quote}
    
     \textbf{Transfer principle}
    
     \begin{center}
         \emph{``the ideal elements          
           correctly determine the real things''}
     \end{center}
  
  \begin{quote}
     Perhaps \emph{the infinite and infinitely small} [numbers] 
    that we conceive \emph{are imaginary, nevertheless} [they are] \emph{suitable 
    to determine the real things}, as usually do the imaginary roots.
    They are situated in the ideal regions, from  where things are ruled by
    laws, even though they do not lie in the part of matter.
    
    (\emph{Letter to Johann Bernoulli}, 1698)

     %
    
\end{quote}

  In this paper, we try and give a precise mathematical formulation 
  of these principles in the context of
  the \emph{Topological Extensions} of \cite{DF05}, structures     
   which generalize at once compactifications, completions, and nonstandard          %
   models (see also \cite{BDF06}). 
   
   Given a set $M$, a 
   \emph{topological extension} of $M$ is a $T_1$ space $\hM$, where 
   $M$ is a dense 
   subspace and every function $f:M\to M$ has a distingushed continuous 
   extension $\f: \hM\to\hM$ that preserves compositions and local 
   identities. The operator $\ast$ can be appropriately defined so as to 
   provide also
   all properties $P$ and relations  $R$ with  distinguished 
   extensions $\hP,\hRR$ to $\hM$.
   
  Following the basic idea that the elements of the [``standard''] set $M$ are 
  the ``real objects'' 
  of the ``actual world'', whereas the extension $\hM$ 
  contains also the ``ideal elements'' of all ``possible worlds'', 
  an appropriate interpretation of 
  the Leibniz's principles   in the context of topological extensions
  might be
  
  \begin{description}
     \item[Ind] \emph{different elements of $\,\hM$ are separated by the 
     extension $\hP$ of some property $P$ of $M$;} 
 
     \item[Poss] \emph{if  the extensions 
     of a family $\F$ of properties of $M$ are not simultaneously 
     satisfied in $\hM$, then there are finitely 
     many properties of $\F$ that are not simultaneously satisfied in 
     $M$;}
  
     \item[Tran] \emph{a statement involving elements, properties 
     and relations of $M$ is true  if and only if the corresponding 
     statement about their extensions is true in 
     $\hM$.}\footnote{~Clearly one has to admit only ``first order''
     statements, so as to avoid trivial inconsistencies.}
 \end{description}
  
  We shall see that the extended properties correspond exactly to the 
  clopen subsets of $\hM$, and so the above principles turn out to be respectively 
    a
   property of \emph{separation}, of \emph{compactness}, and of 
   \emph{analyticity} of the topology of $\hM$.
Grounding on results of \cite{DF05,BDF06}, we obtain the somehow 
surprising consequence that the Leibnz's principles 
\emph{can be fulfilled in pairs, but not all          
   three together}.

\smallskip
The paper is organized as follows.
In Section \ref{top}, we give the precise definition of \emph{topological
extension} and we recall the main  properties stated in \cite{DF05}.  
In particular, in Subsection \ref{stop}, we introduce the \emph{$S$-topology} 
and we determine its connection with the principles \textbf{Ind} and 
\textbf{Poss}. In
Subsection \ref{can} we study the \emph{canonical map} from a topological
extension of $X$ into
the Stone-$\check{\rm C}$ech compactification $\beta X$ of the 
discrete space $X$: we obtain \emph{inter alia} that  the Stone-$\check{\rm C}$ech
compactification itself is essentally  the \emph{unique} topological 
extension that satisfies \emph{both principles} \textbf{Ind} and  
\textbf{Poss}. 

In Section \ref{hyp}, we present two simple properties that characterize 
all those topological 
extensions (\emph{hyperextensions}) that satisfy the \trans\ \tran. 
A complete characterization of the hyperextensions 
satisfying also \ind\ is derived at once, and with it the 
\emph{impossibility of satisfying simultaneously the three Leibniz's 
principles}. In Subsection \ref{star} we show how to topologize 
arbitrary nonstandard models, so as to obtain also topological 
extensions where the principles \poss\ and \tran\ hold together.

A few concluding remarks and open questions, in particular the set-theoretic 
problems 
originated by the combination of \ind\ with \tran, can be found in the final 
Section \ref{froq}.

\smallskip
In general, we refer to \cite{Eng} for all the topological notions and
facts used in this paper, and to \cite{CK} for definitions and facts
concerning ultrapowers, ultrafilters, and nonstandard models. General 
references for nonstandard Analysis could be \cite{Ke76,nsa}; 
specifical for our ``elementary'' approach is  \cite{BDF06}.

\smallskip
The author is grateful to  Vieri Benci, Mauro Di Nasso and Massimo 
Mugnai for useful discussions
and suggestions.

\section{Topological extensions and the Identity of Indiscernibles}
\label{top}
In this section we review the main features of the \emph{topological extensions}  
introduced in the paper \cite{DF05}; these structures  naturally 
accomodate, within a general unified framework, both  \emph{Stone-$\check{\rm C}$ech  compactifications of 
discrete spaces}  
and  \emph{nonstandard models} (see also \cite{BDF06}). 
The most important characteristic shared by compactifications 
and completions in topology, 
and by nonstandard models of analysis is  the existence of a 
distinguished extension $\f:\X\to \X$ for each 
function $f:X\to X$.  Given an arbitrary  set $X$, we consider here a 
topological extension of $X$ as a sort of ``topological completion'' 
$\,\X$, where
  the ``$\ast$'' operator provides a 
\emph{distinguished continuous extension} of each  function $f:X\to X$.

\begin{Definition}\label{text}
\emph{The $T_{1}$ topological
space $\X$ is a
\emph{topological extension} of $X$ if $X$ is a \emph{discrete dense} subspace of
$\X$, and  a 
\emph{distinguished continuous extension} $\f:\X \to \X$ is associated to
each function $f:X \to X$, so as to satisfy the following conditions:}

\begin{itemize}
\item[$(\mathsf{c})$]  
$\g\circ \f = \ns(g\circ f)$~  for all
$f,g:X \to X$, \emph{and}
\item[$(\mathsf{i})$] 
if $f(x)=x$ for all $x\in A\incl X$, then
$\f(\xi)=\xi$ for all $\xi \in\_A$.
\end{itemize}

\end{Definition}

Since a finite
set cannot have \emph{nontrivial} topological
extensions, we are interested only in infinite sets, and
for convenience we stipulate that
$\N \incl X$.

It is easily seen that the operator $\ast$ preserves also constant and
characteristic functions (see  Lemma 1.2 of \cite{DF05}). So,
by using the characteristic functions, the operator $\ast$ 
provides also  an {\em extension} $\hA$ for every \emph{subset}  
$A\incl X$, which turns out to be a \emph{clopen} superset of $A$, and 
actually the \emph{closure} 
$\_ A$ of $A$ in $\,\X$.

Notice that, if the topological extension $\X$ of $X$ is  Hausdorff,
then $\f$ is the
\emph{unique} continuous extension of $f$, because $X$ is dense. Therefore
 properties
$(\mathsf{c})$ and $(\mathsf{i})$ are automatically satisfied
(see \cite{BDF02}, where Hausdorff topological extensions  have
been introduced and studied).
However considering only Hausdorff
spaces would have turned out too restrictive: we shall see below that 
the
Hausdorff topological extensions of $X$ are particular
subspaces of the Stone-$\check{\rm C}$ech compactification $\beta X$
of the discrete space $X$ that, in general, are not 
nonstandard extensions. In fact, the existence of Hausdorff 
nonstandard extensions, although consistent, has not yet 
been proved in \zfc\ alone
(see \cite{DF05,DF06}) and Section \ref{froq} below.  These are the reasons why
we only require  that topological extensions be $T_1$ spaces.

\medskip
\subsection{The $S$-topology}
\label{stop}

In order to study our versions of the Leibniz's principles for
 topological 
extensions, it is useful to consider on ~$\X$ the so called 
\emph{$S$-topology},\footnote{~The  
$S$-topology (for \emph{Standard} topology)
is  a classical notion of nonstandard
 Analysis, already considered since \cite{rob66}.} \ie\
the topology
    generated by the
(clopen) sets ~$\hA = \_ A$ for $A\incl X$.
The $S$-topology  is
obviously  \emph{coarser than} or \emph{equal to} 
the original topology of $\, \X$, and 
we can characterize the respective \emph{separation properties} as 
in Theorem 1.4 of \cite{DF05}:

\begin{Theorem}\label{sep}
    Let $\X$ be a topological extension of $X$. Then
\begin{enumerate}
    \item  The $S$-topology of $\, \X$ is either $0$-dimensional or
not $T_{0}$.

\item  $\X$ is Hausdorff if and only if the  $S$-topology is $T_{1}$,
hence $0$-dimensional.

\item  $\X$ is regular if and only if the  $S$-topology
is the topology of $\, \X$ (and so $\,\X$ is $0$-dimensional).
\end{enumerate}
\end{Theorem}

{\bf Proof.}~

\noindent$1.$~ The $S$-topology has a clopen basis by definition. In
this topology the closure of a point $\, \xi$ is $\, M_{\xi} = \bigcap
_{\xi \in \_ A} \_ A$.
If  $M_{\xi} = \{\xi \}$ for all $\,\xi\in\X$, then the $S$-topology is $T_{1}$,
hence $0$-dimensional.
Otherwise let $\eta \ne \xi$ be in  $M_{\xi}$.
Then $\eta$ belongs to the same clopen sets as $\xi$, and the
$S$-topology is not $T_{0}$. In fact, given
$A\incl X$, $\xi \in \_A $ implies
         $\eta \in \_A $, by the
choice of $\eta$. Similarly  $\xi \notin \_A $ implies
$\xi \in\! \_{X\/ A} $, hence  $\eta \in\! \_{X\/ A} $ and
$\eta \notin \_A $.

\smallskip
\noindent$2.$~ By point 1, the $S$-topology is Hausdorff (in fact
$0$-dimensional)
    whenever it is $T_{1}$. Therefore also the topology of $\X$  is Hausdorff,
being finer than the $S$-topology.
For the converse, let $U,V$ be disjoint neighborhoods of the points
$\xi, \eta \in\X$,  and put
$A=U\cap X$, $B=V\cap X$. Then $\xi \in \_A$,
$\eta \in \_B$, and  $\_B \cap \_A =\0$. Therefore $\eta \notin
M_{\xi}$, and the $S$-topology is $T_{1}$.

\smallskip
\noindent$3.$~ The closure of an open subset $U\incl \X$ is the clopen
set $\_{U\cap X}$. Therefore any closed neighborhood  of $\xi\in
\X$  includes a clopen one. Since the clopen sets are a basis of the
$S$-topology, $\, \X$ can be regular if and only if its original 
topology is  the $S$-topology
(and so the latter is $T_{1}$, hence $0$-dimensional).

\hfill $\Box$

\bigskip
Now the principle \textbf{Ind} simply means 
that  the $S$-topology of $\X$ is Hausdorff. On the other hand, the principle 
\textbf{Poss} states that every proper filter of clopen sets 
has nonempty intersection, \ie\ that the $S$-topology of $\X$ is 
quasi-compact.\footnote{~Following \cite{Eng}, we call compact only 
Hausdorff spaces.}
So we have

\begin{corollary}\label{lp12}
 Let $\X$ be a topological extension of $X$. Then
 \begin{enumerate}
     \item  the principle \emph{\textbf{Ind}} holds if and only if $\,\X$ is 
     Hausdorff;
 
     \item  the principle \emph{\textbf{Poss}} holds if and only if the 
     $S$-topology of $\,\X$ is 
quasi-compact;
 
     \item  both principles \emph{\textbf{Ind}} and  \emph{\textbf{Poss}} 
     hold in 
     $\X$ if and only if the 
     $S$-topology of $\X$ is compact. So either $\X$ is compact, or 
     it is not regular, but becomes compact by suitably weakening 
     its topology, still maintaining all functions $\f$ 
     continuous.
 \end{enumerate}
 \end{corollary}

{\bf Proof.}~
We only need to prove the last assertion of Point 3. If $\X$ is regular, then 
its topology agrees with the $S$-topology. If not, then all functions 
$\f$ are continuous also with respect to the coarser $S$-topology, 
because the inverse images of clopen sets are clopen.

\hfill $\Box$

\medskip
\subsection{The canonical map and the principle \emph{\textbf{Ind}}}
\label{can}

Any topological extension
of $X$ is canonically mappable into the \emph{Stone-$\check{\rm C}$ech
compactification} $\beta X$ of the discrete space $X$.\footnote{~For
various definitions
and properties of the Stone-$\check{\rm C}$ech
compactification see \cite{Eng}.} 

If $X$ is a
discrete space, identify $\beta X$  with the set of all
ultrafilters over $X$, endowed with the topology having as basis
$\{ {\mathcal O}_{A} \mid A\incl X \}$, where ${\mathcal O}_{A} $ is the
set of all ultrafilters containing A. So the embedding  $\, e:X \to \beta
X$ is given by the principal ultrafilter $$e(x)=\{A\incl X\mid x\in 
A\},$$ and 
the unique continuous extension $\overline f:\bg X\to \bg X$ of $f:X\to X$
 can be defined by putting $$ \_f(\U)=\{A\incl X\mid
 f^{-1}(A)\in \U\}.$$

Given a topological
extension $\X$ of $X$ and a point $\xi\in\X$, put  
$$\U_{\xi} = \{ A\incl X \mid \xi \in \hA\},$$ which 
is an
ultrafilter over $X$, and define the \emph{canonical  map}
$$\upsilon : \X \to \beta X\ \
\mbox{by}\ \ \upsilon (\xi) =
\U_{\xi}.$$
Then we can reformulate  Theorem 2.1 of \cite{DF05} in terms of the 
Leibniz's principles \textbf{Ind} and \textbf{Poss}. Namely

\begin{Theorem}\label{tcan}~Let 
    $\X$ be a topological extension of $X$, and let $\bg X$ be the 
    Stone-$\check{\rm C}$ech
compactification of $X$. Then
\begin{enumerate}
\item The canonical map
$\upsilon : \X \to \beta X$  is the unique continuous
extension to $\X$ of the embedding  $\, e: X \to \beta X$, and
$$\upsilon\circ\f = \_f\circ\upsilon\ \  \mbox{for all}\  f:X\to
X.$$

\item
The map $\upsilon$ is injective if and only if \emph{\textbf{Ind}} holds in 
$\X$.
\item The map $\upsilon$ is surjective if and only if
\emph{\textbf{Poss}} holds in  $\,\X$.

\end{enumerate}

\end{Theorem} 

{\bf Proof.}~

\noindent 1.~For all $x\in X$, $\U_{x}$ is the principal ultrafilter
generated by $x$, hence $\upsilon$ induces the canonical
embedding of  $X$ into $\beta X$. If ${\cal O}_{A} $ is a basic open
set of $\beta
X$, then
$\upsilon^{-1}({\cal O}_{A}) = \_A$, hence $\upsilon $ is continuous
w.r.t.~the $S$-topology, and \emph{a fortiori} w.r.t.~the (not coarser)
topology of $\X$.
On the other hand, let a continuous map  $\varphi : \X \to \beta X$ be
given. Since
${\cal O}_{A} $ is clopen, also
	 $\varphi^{-1}({\cal O}_{A}) $ is clopen and so it is
the closure
$\_B$ of some
	 $B\incl X$. If $\varphi$ is the identity on $X$, then $\_B \cap X =A$,
	 hence $B=A$ (see Lemma 1.2 of \cite{DF05}). Therefore all points of
$M_{\xi}= \{\eta\in\X\mid \forall A\incl 
X\,(\xi\in\hA\,\Imp\,\eta\in\hA)\,\}$
	 are mapped by $\varphi$ onto $\upsilon(\xi)$, 
	 and so $\upsilon = \varphi$.

     Moreover, for all $\, \xi \in \X$, one has $\,\xi \in \_A \ \Leftrightarrow \
	 \f(\xi) \in \_{f(A)}$,   or equivalently
	 $A\in {\cal U}_{\xi} \ \Leftrightarrow \ f(A) \in {\cal U}_{\f(\xi)}$
	 (see  Lemma 1.3 of \cite{DF05}); hence
	  $\overline f \circ \upsilon = \upsilon \circ \f$,
	 and Poiny 1. is completely proved.

\noindent 2.	 The map $\upsilon$ is injective if and only if the
      $S$-topology is $T_{1}$, and this fact is equivalent
	 to $\X$ being Hausdorff, by Theorem \ref{sep}, or to 
	 \textbf{Ind}, by Corollary 1.3. Moreover in this
	 case $\upsilon$ is a homeomorphism w.r.t.~the $S$-topology, which
	 is the same as the topology of $\X$ if and only if the latter
	 is regular (hence $0$-dimensional).

	 \noindent 3.
	 The map $\upsilon$ is surjective if and only if every maximal filter
	 in the field of all clopen sets of $\,\X$ has nonempty intersection. This is
	 equivalent to every proper filter  having nonempty
	 intersection, which in turn is equivalent to the $S$-topology of $\X$ 
	 being
quasi-compact, i.e.~to the principle \textbf{Poss}, by Corollary 1.3.

\hfill $\Box$

\medskip
Notice that the map $\upsilon$ induces a bijection between the basic
open sets ${\cal O}_{A} $ of $\beta X$ and the clopen subsets $\hA$
of $\,\X$. Therefore $\upsilon$ is open if and only if $\,\X$ has the
$S$-topology.

 Call \emph{invariant} a subspace $Y$ of $\,\X$ (respectively of $\beta X$)
   if  
$$\f(\xi)\in Y\ \mbox{(resp. $\_f(\xi)\in Y$) for all}\ f:X\to X\ 
\mbox{and all} \
\xi\in Y.$$

It is easily seen that any invariant subspace $Y$ of $\,\X$ is itself a
topological extension of $X$, and it is mapped by $\upsilon$ onto an
invariant  subspace of $\beta X$.
If $\, \X$ is homeomorphic to a subspace of $\beta X$, then it is
 $0$-dimensional, hence it has the  $S$-topology, by Theorem \ref{sep}.
         Conversely, if $\, \X$ has the $S$-topology, then $\upsilon$ is
         injective. Moreover, for all $A\incl X$,
$\upsilon (A) =  {\cal O}_{A}\cap \upsilon (\X)$, hence $\upsilon$
is a homeomorphism between $\X$ and its image. On the other hand,
if $\, \X$ is Hausdorff but not regular, then $\upsilon$  is
injective and continuous, but not open.

Whenever $\X$ verifies the principle \textbf{Ind}, the map $\upsilon$ 
can always be turned
into a homeomorphism,
either by endowing $\upsilon (\X)$ with a
suitably finer topology, or by
taking on $\X$ the (coarser) $S$-topology.
So any such extension makes use of \emph{the
same ``function-extending mechanism''} as the Stone-$\check{\rm C}$ech
       compactification. Moreover, if also \textbf{Poss} holds, then 
       $\X$ can be taken to be $\bg X$ itself, possibly endowed with 
       a suitably finer topology.
       
More precisely, the above discussion provides the same characterization 
of all topological extensions  
 satisfying the principle \textbf{Ind} that has been given
in  Corollary 2.2 of \cite{DF05}, namely:

\begin{corollary}\label{ind}
 A topological 
   extension $\,\X$ of $X$ satisfies \emph{\textbf{Ind}} if and
only if  the canonical map $\upsilon$
provides a continuous bijection  $($a homeomorphism when $\X$ is 
regular$)$  onto
an invariant subspace of $\, \beta X$.

Moreover $\,\X$ satisfies also \emph{\textbf{Poss}} if and only if 
$\upsilon$ is onto $\bg X$.
\hfill $\Box$
\end{corollary}

\section{Topological hyperextensions and the Transfer 
Principle}\label{hyp}       
       
The Transfer Principle  \tran\ is the very ground of the usefulness of the 
nonstandard methods in mathematics. It allows for obtaining correct 
results about, say, the real numbers by using ideal elements like 
actual infinitesimal or infinite numbers.

In fact, both properties $(\mathsf{c})$ and $(\mathsf{i})$ of 
Definition \ref{text} are 
instances of the \trans, for they correspond  to the 
statements
$$\forall x\in X\, .\, f(g(x))=(f\circ g)(x)\ \
\mbox{and}\ \
\forall x\in A\, .\, f(x)=x,$$
respectively. So all topological extensions already satisfy several 
important cases of the
\trans.
\Pes, if $f$ is constant, or injective, or surjective, then
so is  $\f$. More important, we have already used the fact that
the extension of the characteristic function of
any subset $A\incl X$ is the characteristic function of the closure 
$\_ A$ of
$A$ in $\X$, thus we can put $\hA =\_ A$ and obtain a 
Boolean isomorphism between the field $\P(X)$ of all subsets of $X$ 
and the field $\C\ell(\X)$ of all clopen subsets of $\X$ (see \cite{DF05},
Lemmata 1.2 and 1.3).     

On the other hand, many basic cases of  the \trans\ may fail, 
because topological extensions comprehend, besides 
nonstandard models, also all invariant subspaces of the Stone-$\check{\rm C}$ech
       compactifications of discrete spaces.
In order to obtain the full principle \tran,
we postulated in \cite{DF05} two
additional properties, namely

\begin{Definition}\label{texan}
The topological extension $\X$ of $X$ is a
\emph{hyperextension}\footnote
{~Topological hyperextensions are
in fact \emph{hyper-extensions} in the sense of 
\cite{BDF06} (\ie\ nonstandard models), by Theorem \ref{eltop} below.} 
if
\begin{itemize}
\item[$(\mathsf{a})$] ~for all $f,g:X\to X$\\
~~~~~~$f(x) \ne g(x)$ for all $x\in X$
$\ \ \Longleftrightarrow \ \ \, \f(\xi) \ne \g(\xi)$ for all $\xi\in \X$;

\smallskip
\item [$(\mathsf{p})$]~there
exist $p,q:X\to X$ such that\\
~~~~~~for all $\xi,\eta \in \X$ there
exists $\zeta \in \X$ such that
$\ \xi = \p(\zeta)$ and  $\eta = \q (\zeta)$.
\end{itemize}

\end{Definition}

The property $\mathsf{(a)}$, called
\textit{analyticity} in \cite{DF05}, isolates a
fundamental feature that marks the difference between
\emph{nonstandard} and \emph{ordinary continuous} extensions of functions:
\emph{``disjoint functions have disjoint extensions''}.
It is obtained by \tran\ from the statement 
$\forall x\in X\, .\, f(x) \ne g(x)$, and it can be viewed as
the \emph{empty set case}  
of a general \emph{``principle of preservation of
equalizers''}:
\begin{description}
\item[$(\es)$] {}~~~~~~~~~~~$\{ \xi \in \st X \mid \st f(\xi)= \st g(\xi)\} =
\st{\{ x \in X \mid f(x)= g(x)\}}.$
\end{description}

The  property $\mathsf{(p)}$, called
\textit{coherence} in \cite{DF05}, provides a sort of ``internal coding of
pairs'', useful
for extending  multivariate functions ``parametrically'': this
possibility is essential in order to get the full principle  \tran, which 
involves relations of any arities.\footnote{~
The  \emph{ratio} of considering only unary functions lies in the
following facts that hold in every   topological
hyperextension $\X$ of $X$ (see Section 5 of \cite{DF05}):\\
\emph{- For all $\xi_{1},\ldots,\xi_{n}\in \X$ there exist
$ p_{1},\ldots,p_{n}: X \to X$ and  $\zeta\in \X$ such that
$\p_{i}(\zeta) = \xi_{i}$.}\\
\emph{- If $p_{1},\ldots,p_{n},q_{1},\ldots, q_{n}: X \to X$ and
$\xi,\eta \in \X$ satisfy  $\p_{i}(\xi) = \q_{i}(\eta)$, then\\ 
${}~~~~~~~~~~~~~~~~\ns (F \circ (p_{1},\ldots,p_{n}))(\xi) =
\ns (F \circ (q_{1},\ldots, q_{n}))(\eta)$ for all $F: X^{n} \to X$}.

It follows that there is a unique way of assigning an
extension $\ns F$ to every function $F:X^{n} \to X$ in such
a way that all compositions are preserved.
By using the characteristic functions in $n$ variables one can assign
an extension $\ns R$ also to all $n$-ary relations $R$ on $X$.}
~Notice that  the property  $(\mathsf{p})$ could seem \emph{prima
facie} an illegal instance 
of the \transfer, for it is given in a second order formulation. 
On the contrary, a
\emph{strong uniform version} of that property can be obtained by
fixing $p,q$ as the compositions of a given bijection  $\delta :X\to X\times X \ $
with the ordinary projections $\pi_{1},\pi_{2}: X\times X \to X$,
and then applying \tran\ to the statement
$$\forall x,y\in X\, .\,  \exists z
\in X\, .\ p(z) =x,\ q(z) =y.$$

We shall see below that there are invariant 
subspaces of the Stone-$\check{\rm C}$ech compactification $\bg X$
 where $\mathsf{(a)}$ holds whereas 
$\mathsf{(p)}$ fails and \emph{vice versa}, as well as  
invariant 
subspaces where both fail or hold. So the properties $(\as)$ and $(\ps)$ 
are independent, 
also when \ind\ holds.

We consider very remarkable
 the fact that the combination of \emph{four natural, simple instances} 
 of the \trans, like $(\cs),(\is),(\as)$, and $(\mathsf{p})$,
gives to topological  hyperextensions the
strongest \transfer\ \tran. 
In reason of its importance, we have already given 
three different proofs of this fact in preceding papers of ours: a 
``logical'' and a ``logico-algebraic'' 
proof in  \cite{DF05}, and 
a ``purely algebraic'' proof in \cite{F03} (see also the survey in 
\cite{BDF06}). So we state here without proof
the following theorem: 

\begin{Theorem}[]\label{eltop}~ 
A topological extension $\X$ of $X$ satisfies the principle \emph{\tran}\ if 
and only if it is a hyperextension.
\qed
\end{Theorem}

We are now able to characterize all topological
extensions satisfying both principles \textbf{Ind} and \tran. 
These extensions are spaces of ultrafilters, according to Corollary 
\ref{ind}. So we  use the reformulation  in terms of ultrafilters
given in \cite{DF05} for the condition $(\mathsf{e})$ above.

Call an ultrafilter $\U$ on $X$  \emph{Hausdorff}\footnote{~The
      property $(\mathsf{H})$ has been introduced in \cite{dt} under
      the name $(C)$. Hausdorff ultrafilters are studied in
      \cite{DF06} and \cite{bs}.}
      if, for all $f,g:X\to X$,

\smallskip
\noindent $(\mathsf{H}) {}~~~~~~ $
     $ \ \ \,\_f(\U)= \_g(\U)\ \Longleftrightarrow \
      \{\,x\in X \mid f(x)=g(x)\} \in \U .$

\bigskip
Call \emph{directed} a subspace $Y$  of $\beta X$ where the
property $(\ps)$ holds, \ie\ there
exist  $p,q:X\to X$ such that {for all $\U,\V \in Y$ there
exists $\W \in Y$  such that
$\ \U = \_ p(\W)$ and  $\V = \_ q (\W)$}.

By combining Theorem \ref{eltop} with Corollary \ref{ind} we obtain
\begin{Theorem}\label{lp13}
  A  topological
      extension $\X$ of $X$ satisfies both principles
      \emph{\textbf{Ind}} and \emph{\tran} if and
only if  the canonical map $\upsilon$
is a continuous bijection between $\X$ and
a directed invariant  subspace of $\beta X$ that contains only 
Hausdorff ultrafilters. 

Moreover $\yg$ is a homeomorphism
if and only if   $\X$ has the $S$-topology, or equivalently is regular.

\qed
\end{Theorem}

Now it is easy to show that the properties $(\as)$ and $(\ps)$ are 
independent.

For $\U\in\bg X$ let
$Y_{\U}=\{\_f(\U)\mid f:X\to X\,\}$ be the invariant subspace 
generated by $\U$. Clearly $Y_{\U}$ is directed, so $(\ps)$ holds for 
all ultrafilters $\U$, whereas $(\as)$ holds if and only if $\U$ is 
Hausdorff.
On the other hand, let $\U$ and $\V$ be Hausdorff ultrafilters such 
that  neither 
of them belongs to the invariant subspace generated by the other one: 
then $Y_{\U}\cup Y_{\V}$ is an invariant subspace where $(\as)$ holds,
but it is not directed, hence $(\ps)$ fails.

We shall deal in the final section with the set theoretic strength of 
the combination of \ind\ with \tran. By now we simply recall that 
there are plenty of non-Hausdorff ultrafilters (\pes\ all diagonal tensor 
products $\U\otimes\U$). Thus we can easily conclude
\begin{corollary}\label{lp123}
 No topological 
   extension satisfies at once the three Leibniz's principles
   \emph{\textbf{Ind}}, \emph{\textbf{Poss}}, and \emph{\tran}.

\hfill $\Box$
\end{corollary}

\subsection{The star topology}\label{star}

We are left with the task of combining \poss\ with \tran.
To this aim we recall that a nonstandard model whose $S$-topology is 
quasi-compact is commonly called \emph{enlargement}. It is well known 
that every structure has arbitrarily saturated 
 elementary extensions (see \pes\ \cite{CK}), and obviously a 
 $2^{|X|^{+}}$-saturated extension of $X$ is an enlargement (see \pes\
 \cite{nsa} or \cite{BDF06}). 
 So, if we can topologize every nonstandard extension of $X$ in such 
 a way that all functions $\f$ become continuous, then we get a lot 
 of topological hyperextensions satisfying \poss.
 
 This task has been already accomplished in \cite{DF05}, where the 
  the \emph{coarsest} such topology is defined.

Every topological extension $\X$ should be a $T_{1}$ space, so 
all sets  of the form
$E(f,\eta) =\{ \xi \in \X \mid \f(\xi) = \eta\}$, for
$f:X\to X$ and $\eta
    \in \X$,
should be closed in $\X$.
The (arbitrary) intersections of finite unions of such sets are the
    \emph{closed} sets of a topology, called the \emph{Star topology},
     which  is by construction the \emph{coarsest $T_{1}$
    topology on $\X$ that makes all
    functions $\f$ continuous}.

When $\X$ is a nonstandard extension of $X$, the four defining properties 
$(\cs),(\is),(\as)$, and $(\ps)$ of  topological hyperextensions 
are fulfilled by hypothesis. So one has only to prove that $X$ is
dense in  $\X$ in order to obtain
\begin{Theorem}[Theorem 3.2 of \cite{DF05}]\label{start}~ 
Any nonstandard extension $\,\X$ of $X$, when  
equipped with the Star topology,
becomes a topological hyperextension of $X$. 
Conversely, any topological hyperextension $\X$ of $X$ 
is a nonstandard extension, possibly endowed 
with a topology finer than the Star topology.
\end{Theorem}

{\bf Proof.}~
We have only to prove that $X$ is dense.

If $X\incl \bigcup_{1\le i\le
n}E(f_{i},\eta_{i})$, then
 we may consider w.l.o.g. only those
components with $\eta_{i}\in X$,
because
each  $\f_{i}$ maps any point $x\in X$
to the point $f_{i}(x)\in X$. Hence
the \transfer\ of the nonstandard extensions may be applied to the 
statement
 $$\forall x\in X \ (f_{1}(x)=\eta_{1} \vee \ldots
\vee f_{n}(x)=\eta_{n}),$$
thus producing
$$\forall \xi\in \X \
(\f_{1}(\xi)=\eta_{1} \vee \ldots
\vee \f_{n}(\xi)=\eta_{n}).$$
So the whole space $\,\X$ is included in $\bigcup_{1\le i\le
n}E(f_{i},\eta_{i})$, and $X$ is dense.

\qed

\bigskip
So, in order to get topological extensions satisfying both principles 
\poss\ and \tran, we have only to put the \emph{star topology} on any 
\emph{nonstandard enlargement} of $X$.

\section{Final remarks and open questions}\label{froq}

  We have seen that (at least) one of the three principles that we 
  have investigated has to be left out. The most reasonable choice seems 
  to be that of dropping \ind. In fact, even if one neglects the set 
  theoretic problems that will be outlined below, one should pay 
  attention to
  Leibniz himself.
  
   \begin{quote}

[...] cette supposition de deux indiscernables [...] paroist possible 
en termes abstraits, mais elle n'est point compatible avec l'ordre 
des choses [\ldots]

Quand je nie qu'il y ait [\ldots] deux corps indiscernables, je ne 
dis point qu'il soit impossible 
absolument d'en poser, mais que c'est une chose contraire ˆ la 
sagesse divine [\ldots]

Les parties du temps ou du lieu [...] sont des choses ideales, ainsi
elles se rassemblent parfaitement comme deux unit\'es abstraites. 
Mais il n'est pas de m\^eme de deux Uns concrets [...] c'est \`a 
dire veritablement actuels.

Je ne dis pas que deux points de l'Espace sont un meme point, ny 
que deux instans du temps sont un meme instant comme il semble
qu'on m'impute [\ldots]\footnote{\ldots this supposition of two indiscernibles \ldots seems 
abstractly possible, but it is incompatible with the order of 
things\ldots\\
When I deny that there are  \ldots two indiscernible bodies, I do not 
say that [this existence] is absolutely impossible to assume, but 
that it is a thing contrary to Divine Wisdom \ldots \\
The parts of time or place \ldots are ideal things, so 
they perfectly resemble like
two abstract unities. But it is not so with two concrete Ones,\ldots
that is truly actual [things].\\
I don't say that two points of Space are one same point, neither 
that two instants of time are one same instant as it seems that one 
imputes to me \ldots}

(\emph{Fifth letter to Clarke}, \cite{LC}, pp. 131-135)
   \end{quote}

   It appears that  Leibniz  considered the identity of indiscernibles 
   as a ``physical'' rather than a ``logical'' principle: it is 
   actually true, but its negation is non-contradictory in principle,
   so \emph{it could fail in some 
   possible world}. 
   Moreover only ``properties of the real world'' $M$
   are considered in all these principle: so it seems natural, and 
   not absurd, to assume that objects indiscernible by these ``real'' 
   properties may be separated by some abstract, ``ideal'' property 
   of $\hM$.
   
On this ground we finally decide  to call \emph{Leibnizian} a topological 
extension that satisfies both \textbf{Poss} and \textbf{Tran}, and so
necessarily not \textbf{Ind}. Thus the existence of plenty of 
\emph{Leibnizian extensions} is granted by the final results of Section 
\ref{hyp}, without any need of supplementary set theoretic hypotheses.

\subsection{Existence of Hausdorff  extensions}\label{shaus}

As shown by Theorem \ref{lp13}, combining \ind\ with \tran\
requires  special ultrafilters, named \emph{Hausdorff} in 
Section \ref{hyp}. Despite the apparent weakeness of their defining 
property $(\mathsf{H})$, which is actually true whenever any of
the involved functions is injective (or constant), not much is known about 
Hausdorff ultrafilters.

On  \emph{countable} sets, the property
$(\mathsf{H})$ is satisfied by \emph{selective}
ultrafilters as well as by \emph{products of pairwise nonisomorphic 
selective}
ultrafilters (see \cite{DF06}), but their existence in pure \zfc\ is still unproved.
However any hypothesis providing infinitely many 
nonisomorphic 
selective ultrafilters over $\N$, like the Continuum Hypothesis 
$\mathsf{CH}$ or Martin Axiom $\mathsf{MA}$,
provides infinitely many non-isomorphic hyperextensions of $\N$ that 
satisfy \ind.

On \emph{uncountable} sets the situation is highly problematic: 
it is proved in \cite{DF06} that  Hausdorff ultrafilters on sets of size
not less than $\uf$ cannot be \emph{regular}.\footnote{
~~$\uf$ is the least size of an ultrafilter basis on $\N$. All what is provable in 
$\mathsf{ZFC}$ about the size of $\uf$
is that $\aleph_{1} \le \uf \le 2^{\aleph_{0}}$ 
(see e.g. \cite{bl03}).} 
In particular, the existence of a hyperextension 
satisfying \ind\ with uniform 
ultrafilters, even on $\R$, would  imply that of inner models 
with measurable 
 cardinals. (To be sure, such ultrafilters have been obtained only by 
 much stronger hypotheses, see \cite{kt}).

Be it as it may, as far as we do not abide $\mathsf{ZFC}$ as our foundational
theory, \emph{we cannot prove that
hyperextensions without indiscernibles exist at all}.
 
\subsection{Some open questions}\label{sopq}

We conclude this paper with a few open questions that involve special
ultrafilters, and so should
 be of independent set theoretic interest.

\smallskip
\begin{enumerate}
    \item Is the existence of topological hyperextensions of $\N$ without 
indiscernibles
provable in $\mathsf{ZFC}$, or at least 
derivable from
set-theoretic hypotheses weaker than those providing selective 
ultrafilters? \Pes\ from ${\mathfrak x = \mathfrak c}$, where  
$\mathfrak x$ is a
 cardinal 
invariant of the continuum not dominated by \textbf{cov}$(\B)$?

    \item  Is it consistent with $\mathsf{ZFC}$ that there are 
    nonstandard real lines 
$\hR$ without indiscernibles where all ultrafilters 
are uniform?

    \item  Is the existence of countably compact hyperextensions
consistent with $\mathsf{ZFC}$? (These extensions would be of great 
interest, because they would verify \ind, 
\tran, and the weakened version of \poss\ that considers only 
\emph{sequences} of properties.)
\end{enumerate}

 %








\end{document}